\pdfoutput=1

\documentclass{article}
\usepackage[14pt]{extsizes}
\usepackage[utf8]{inputenc}
\usepackage[english,russian]{babel}
\usepackage{math}
\usepackage{xcolor}
\usepackage{amssymb}
\usepackage{amsmath}
\usepackage{graphics}
\usepackage{tikz}
\usepackage{pgfplots}
\usepackage{comment}
\usepackage[ruled, linesnumbered, noend]{algorithm2e}

\title{Внешние биллиарды вне правильного двенадцатиугольника: компьютерное доказательство периодичности почти всех орбит и существования апериодической орбиты}
\author{Филипп Рухович}

\begin{document}

\maketitle

Раздел: теория динамических систем

517.938 УДК

\begin{abstract}
Доказано существование апериодической орбиты для внешнего биллиарда вне правильного двенадцатиугольника, а также что почти все траектории такого внешнего биллиарда являются периодическими; явно выписаны все возможные периоды. Доказательство базируется на фактах, установленных с помощью алгоритма, реализованного автором и выполненного на компьютере. Все вычисления в алгоритме были выполнены абсолютно точно; такие вычисления можно проделать <<вручную>>, однако объем вычислений оказывается для этого слишком большим.
\end{abstract}

\section{Введение}
Пусть $\gamma$ - выпуклая фигура, а $p$ — точка вне ее. Проведем правую относительно $p$ касательную к $\gamma$; определим $Tp \equiv T(p)$ как точку, симметричную $p$ относительно точки касания.

\begin{Def} \label{base}
Отображение $T$ называется внешним биллиардом; фигура $\gamma$ называется столом внешнего биллиарда.
\end{Def}
\begin{Def}
Точку $p$ вне фигуры $\gamma$ назовем периодической, если существует такое натуральное $n$, что $T^np = p$; минимальное такое $n$ назовем периодом точки $p$ и обозначим как $per(p)$. 
\end{Def}
\begin{Def}
Точку $p$ вне фигуры $\gamma$ назовем апериодической, если она --- не периодическая, а ее траектория бесконечна в две стороны. 
\end{Def}
\begin{Def}
Точку $p$ вне фигуры $\gamma$ назовем граничной, если $T^np$ не определено для некоторого $n \in \mathbb{Z}$.
\end{Def}

В данной статье будем полагать, что $\gamma$ --- выпуклый многоугольник.

Внешние биллиарды были введены Бернардом Нойманном в 1950-х годах и стали популярны в 1970-х благодаря Ю.Мозеру \cite{Moser78}. Внешние биллиарды исследовались рядом авторов\ (см. например, \cite{Tab93}, \cite{Schwartz09}, \cite{DF09}, \cite{BC11} а также монографию \cite{Tab05}). Так, Р.Шварц \cite{Schwartz09} показал, что траектория начальной точки может быть неограниченной, тем самым разрешив вопрос Мозера\ -\ Нойманна, поставленный в \cite{Moser78}.

В центре нашего внимания находятся следующие открытые в общем случае \underline{проблемы периодичности}:

\begin{enumerate}
\item Существует ли апериодическая точка для внешнего биллиарда вне правильного $n$-угольника?
\item Какова мера периодических орбит внешнего биллиарда вне правильного $n$-угольника?
\end{enumerate}

Данная статья является продолжением статьи \cite{Rukhovich18}, в которой в деталях исследовался внешний биллиард вне правильного восьмиугольника. С.Л.Табачников в \cite{Tab93} решил проблемы периодичности для случаев $n = 3, 4, 6$: для них апериодической точки нет (а периодические, соответственно, образуют множество полной меры), а также для  $n = 5$ - здесь апериодические точки существуют, но их мера равна нулю. В монографии \cite{Tab05}, опубликованной Американским математическим обществом в 2005 г., С.Л.Табачников приводит результаты компьютерного моделирования для восьмиугольника, но пишет, что <<строгого анализа до сих пор нет>>, и что для других случаев результатов нет. В дальнейшем правильный пятиугольник и связанная с ним символическая динамика подробно исследовались в работе N.Bedaride и J.Cassaigne \cite{BC11} (см. также их монографию \cite{BC11a}).

В монографиях \cite{Schwartz10}, \cite{Schwartz14} Р.Шварц исследовал внешний биллиард вне правильного восьмиугольника и множество связанных с ним вопросов; однако решения проблем периодичности им получено не было.

Правильные 3-х, 4-х и 6-угольники являются простейшими. По мнению сообщества, следующими по сложности с точки зрения проблем периодичности являются случаи $n = 5,10,8,12$, ибо только в этих случаях, по-видимому, существует самоподобие периодических структур. Случай $n = 5$ был исследован Табачниковым \cite{Tab93}; случай $n = 8$ был исследован автором в \cite{Rukhovich18}; случай же $n = 10$ похож на случай $n = 5$ (см. \cite{BC11}). Внешний же биллиард вне правильного двенадцатиугольника обладает наиболее сложной из описанных случаев структурой, и самоподобие здесь достигается более нетривиальным методом, нежели в других случаях. Таким образом, исследование случая $n = 12$ в некотором смысле полностью закрывает проблемы периодичности для <<самоподобных>> случаев.


Основным результатом данной работы являются следующие теоремы.
\begin{Th} \label{MainTreorem}
Для внешнего биллиарда вне правильного двенадцатиугольника существует апериодическая точка.
\end{Th}

\begin{Th} \label{MainTreorem2}
В случае внешнего биллиарда вне правильного двенадцатиугольника, периодические точки образуют вне стола множество полной меры.
\end{Th}

\begin{Th} \label{MainTreorem3}

Введем следующие матрицы:

\[ 
M_{68} := 
\begin{pmatrix}
1 & 0 & 0 & 0 & 0 & 0 & 0 & 0 \\
0 & 1 & 0 & 0 & 0 & 0 & 0 & 0 \\
0 & 0 & 1 & 0 & 0 & 0 & 0 & 0 \\
0 & 0 & 0 & 1 & 8 & 18 & 13 & 24 \\
0 & 0 & 0 & 0 & 2 & 7 & 14 & 29 \\
0 & 0 & 0 & 0 & 0 & 0 & 0 & 0
\end{pmatrix}, M_{66} := 
\begin{pmatrix}
0 & 0 & 0 & 0 & 0 & 0 \\
0 & 0 & 0 & 0 & 0 & 0 \\
0 & 0 & 0 & 0 & 0 & 0 \\
0 & 0 & 0 & 0 & 0 & 0 \\
5 & 4 & 3 & 2 & 1 & 0 \\
1 & 1 & 1 & 1 & 1 & 1
\end{pmatrix},
\]

\[
M_{88} :=
\begin{pmatrix}
2 & 2 & 2 & 2 & 20 & 50 & 26 & 50 \\
2 & 2 & 2 & 2 & 20 & 50 & 26 & 50 \\
4 & 4 & 4 & 4 & 42 & 107 & 74 & 145 \\
2 & 2 & 2 & 2 & 20 & 50 & 48 & 94 \\
0 & 1 & 0 & 0 & 0 & 0 & 0 & 0 \\
1 & 0 & 0 & 0 & 0 & 0 & 0 & 0 \\
0 & 0 & 0 & 1 & 8 & 18 & 13 & 24 \\
0 & 0 & 1 & 0 & 0 & 0 & 0 & 0
\end{pmatrix}.
\] Пусть

$F := \{\begin{pmatrix} 0 \\ 1 \\ 0 \\ 0 \\ 0 \\ 0 \\ 0 \\ 0 \end{pmatrix}, 
              \begin{pmatrix} 0 \\ 0 \\ 1 \\ 0 \\ 0 \\ 0 \\ 0 \\ 0 \end{pmatrix},
              \begin{pmatrix} 9 \\ 2 \\ 5 \\ 2 \\ 0 \\ 0 \\ 0 \\ 0 \end{pmatrix},
              \begin{pmatrix} 6 \\ 4 \\ 10 \\ 4 \\ 0 \\ 0 \\ 0 \\ 0 \end{pmatrix},
              \begin{pmatrix} 0 \\ 0 \\ 2 \\ 3 \\ 0 \\ 0 \\ 1 \\ 0 \end{pmatrix},
              \begin{pmatrix} 24 \\ 24 \\ 120 \\ 102 \\ 0 \\ 0 \\ 18 \\ 0 \end{pmatrix},
              \begin{pmatrix} 48 \\ 48 \\ 156 \\ 108 \\ 0 \\ 0 \\ 24 \\ 0 \end{pmatrix},
              \begin{pmatrix} 4 \\ 4 \\ 9 \\ 4 \\ 0 \\ 0 \\ 1 \\ 0 \end{pmatrix},
              \\
              \begin{pmatrix} 1 \\ 0 \\ 0 \\ 0 \\ 0 \\ 0 \\ 0 \\ 0 \end{pmatrix},
              \begin{pmatrix} 0 \\ 0 \\ 1 \\ 1 \\ 0 \\ 0 \\ 0 \\ 0 \end{pmatrix},
              \begin{pmatrix} 2 \\ 2 \\ 4 \\ 2 \\ 0 \\ 0 \\ 1 \\ 0 \end{pmatrix},
              \begin{pmatrix} 0 \\ 0 \\ 7 \\ 8 \\ 0 \\ 0 \\ 1 \\ 0 \end{pmatrix},
              \begin{pmatrix} 6 \\ 6 \\ 13 \\ 6 \\ 0 \\ 0 \\ 2 \\ 0 \end{pmatrix}
            \},\ G := \{ \begin{pmatrix} 0 \\ 0 \\ 0 \\ 0 \\ 1 \\ 0 \end{pmatrix},
                      \begin{pmatrix} 0 \\ 0 \\ 0 \\ 1 \\ 0 \\ 0 \end{pmatrix},
                      \begin{pmatrix} 0 \\ 0 \\ 0 \\ 24 \\ 36 \\ 0 \end{pmatrix},
                      \begin{pmatrix} 0 \\ 0 \\ 0 \\ 18 \\ 36 \\ 0 \end{pmatrix},
                      \\
                      \begin{pmatrix} 0 \\ 0 \\ 0 \\ 1 \\ 2 \\ 0 \end{pmatrix},
                      \begin{pmatrix} 0 \\ 0 \\ 0 \\ 1 \\ 1 \\ 0 \end{pmatrix},
                      \begin{pmatrix} 0 \\ 0 \\ 0 \\ 2 \\ 1 \\ 0 \end{pmatrix},
                      \begin{pmatrix} 0 \\ 0 \\ 0 \\ 1 \\ 3 \\ 0 \end{pmatrix} \}
            $.
            
Пусть $H := \{M_{66}^kM_{68}M_{88}^nf\ |\ f \in F, k, n \in \mathbb{Z}_{\geq 0} \} \cup
            \{M_{66}^kg\ |\ g \in G, k \in \mathbb{Z}_{\geq 0} \} $, и пусть
            $B := \{12\frac{\begin{pmatrix} 1 & 1 & 1 & 1 & 1 & 1 \end{pmatrix}h}{\text{НОД}(12,\ \begin{pmatrix} 1 & 2 & 3 & 4 & 5 & 6 \end{pmatrix}h)}\ |\ h \in H\}$.
Множество всевозможных периодов точек для внешнего биллиарда вне правильного двенадцатиугольника есть объединение $B \cup \{2*b | b \in B, b\text{ нечетно}\}$.
\end{Th}

В доказательстве всех трех теорем принимают активное участие компьютерные вычисления. Дело в том, что можно ввести систему координат на плоскости таким образом, что вершины двенадцатиугольника и многие связанные с исследованием точки имели бы координаты, лежащие в поле $\mathbb{Q}[\sqrt{3}]$; это, в совокупности с наличием в языке Python типа данных, хранящего целые числа произвольной длины, позволяет проводить необходимые компьютерные вычисления абсолютно точно.   

\section{Базовые обозначения и замечания}

\begin{figure}[h!]
\begin{center}
\includegraphics[width=150mm]{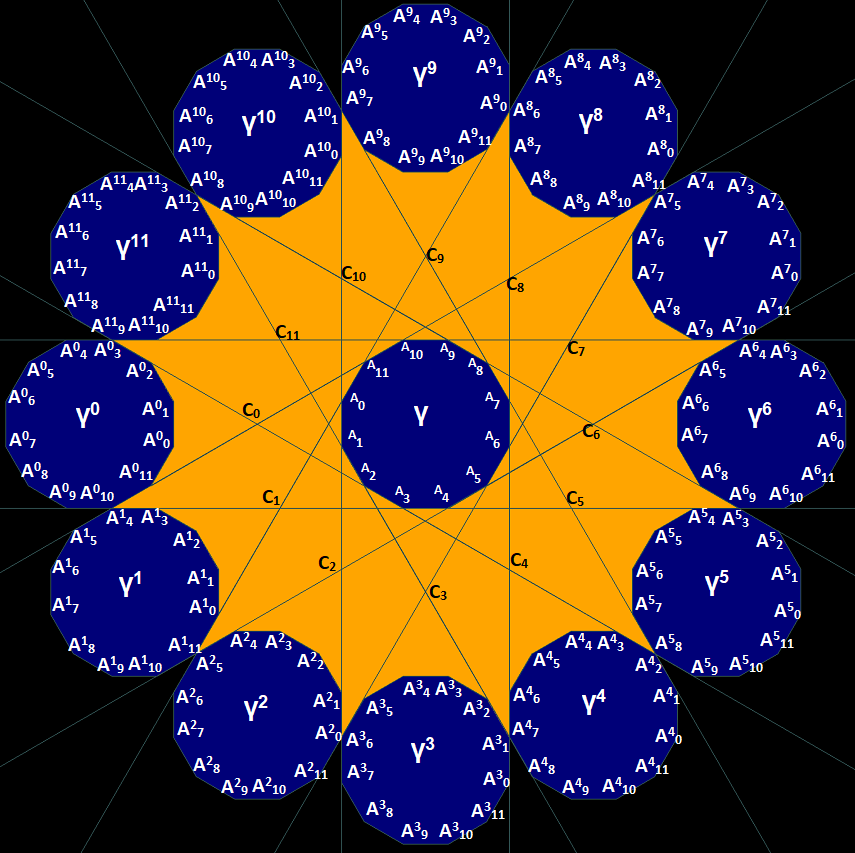}
\caption{Первая инвариантная компонента \cite{Tab93} - фигура $Z$}.
\label{pic:firstInvariantZ}
\end{center}
\end{figure}

Будем следовать плану, намеченному в [Rukhovich18]. Рассмотрим рис. \ref{pic:firstInvariantZ} и изображенный на нем правильный двенадцатиугольник $\gamma = A_0A_1 \ldots A_{11}$. Введем систему координат таким образом, что координаты всех точек лежат в поле $\mathbb{Q}[\sqrt{3}]$. Пусть $l_i$ есть прямая, проходящая через вершины $A_i$ и $A_{(i+1)\ mod\ 12}$, $i = 0, 1, \ldots, 11$, а точка $C_i, i\in[0,11]$ как точку пересечения прямых $l_{(i - 2)\ mod\ 12}$ и $l_{(i+2)\ mod\ 12}$. Пусть $\gamma^i = A^i_0A^i_1\ldots A^i_{11}$ - это открытый многоугольник, симметричный столу $\gamma$ относительно точки $C_i$.

Заметим, что $\forall i \in [0, 12): A^i_1 = A^{(i+1)\ mod\ 12}_6$. Это означает, что многоугольники $\gamma^i$ ограничивают невыпуклый многоугольник

$Z := A^0_3A^0_2A^0_1A^0_0A^0_{11}A^1_4A^1_3A^1_2A^1_1A^1_0\ldots A^{11}_2A^{11}_1A^{11}_0A^{11}_{11}A^{11}_{10}$.

\begin{Lm}
$\forall i \in [0, 12): T(\gamma^i) = \gamma^{(i+5)\ mod\ 12}$.
\end{Lm}

\begin{Lm}
$T(Z) \subset Z \supset T^{-1}(Z)$.
\end{Lm}
 
Также введем <<классическое>> для внешнего биллиарда (см., например, \cite{BC11}) кодирование орбит.

\begin{Def}
Пусть $V_i$, $i \in [0, 12)$, есть угол $A_{(i-1)\ mod\ 12}A_iA^{(i-2)\ mod\ 12}_{i}$.
\end{Def}

\begin{Lm}
Пусть $p$ - произвольная точка вне $\gamma$, а $i \in [0, 12)$. Тогда преобразование $T$ для точки $p$ определено и является центральной симметрией относительно вершины $A_i$, если и только если $p \in \partial V_i$.
\end{Lm}

\begin{Def}
Пусть $p$ - периодическая или апериодическая точка вне стола $\gamma$. Тогда $\rho(p)$ есть бесконечная в обе стороны последовательность $(u_n), n \in \mathbb{Z}$, т.ч. $\forall n \in \mathbb{Z}: T^{n}(p) \in int(V_{u_n})$.
\end{Def}

\begin{Def}
Пусть $p$ - граничная точка вне стола $\gamma$, т.ч. последовательное применение преобразования $T$ может быть выполнено ровно $m \in [0, +\infty]$ раз, а преобразования $T^{-1}$ - $l \in [0, +\infty]$ раз. Тогда $\rho(p)$ есть последовательность $(u_n), n \in [-l, m)$, т.ч. $\forall n \in [-l, m): T^{n}(p) \in int(V_{u_n})$.
\end{Def}

Отметим, что такой код $\rho$ может быть введен для произвольного многоугольного стола. 

\section{Ограничение преобразования}

Заметим, что преобразование $T$ инвариантно относительно поворота на угол $z\frac{\pi}{6}, z \in \mathbb{Z}$ вокруг центра многоугольника $\gamma$. Отождествим точки относительно такого поворота. Как следствие, изучаемую область можно ограничить до угла $A^3_6A_1A_2$, на котором преобразование внешнего биллиарда $T$ индуцирует преобразование $T'$. Фигура же $Z$ ограничивается до фигуры $Z'$, являющейся шестиугольником $A_1A^3_2A^3_3A^3_4A^3_5A^3_6$. Очевидно, $Z'$ похож на <<ракету>> с основанием в многоугольнике $
\gamma^3$; будем называть <<ракетами>> шестиугольники, подобные $Z'$.

\begin{figure}[h!]
\center{\includegraphics[width=150mm]{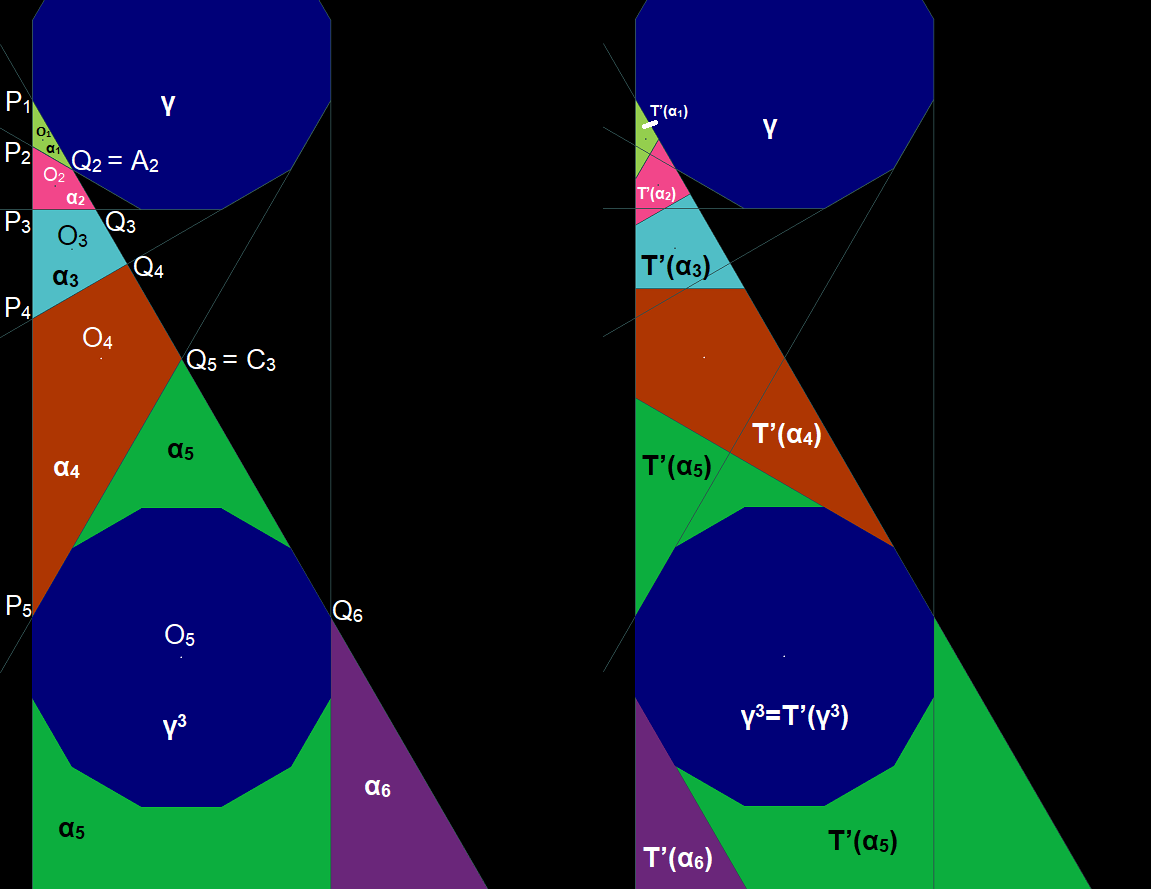} }
\caption{Основные обозначения и преобразование $T'$}.
\label{pic:inducedT}
\end{figure}

Опишем преобразование $T'$ (см. рис. \ref{pic:inducedT}). Пусть точки $P_1$, $P_2$, \ldots, $P_5$ суть точки пересечений луча $A_0A_1$ с лучами $A_2A_1$, $A_3A_2$, \ldots, $A_6A_5$ соответственно, а точки $Q_2$, $Q_3$, \ldots, $Q_6$ суть точки пересечений луча $A_1A_2$ с лучами $A_3A_2$, $A_4A_3$, \ldots, $A_7A_6$ соответственно. Отметим, что некоторые их точек получили вторые, а то и третьи имена; так, $P_1 = A_1$, $Q_2 = A_2$, $Q_5 = C_3$, $P_5 = A^3_6$, $Q_6 = A^3_1 = A^4_6$. Пусть фигура $\alpha_1$ есть треугольник $P_1P_2Q_2$, $\alpha_2$ --- четырехугольник $P_2P_3Q_3Q_2$, $\alpha_3$ --- четырехугольник $P_3P_4Q_4Q_3$, $\alpha_4$ --- четырехугольник $P_4P_5Q_5Q_4$,
$\alpha_5$ --- бесконечная фигура, ограниченная лучом $A^3_6A^3_7 = P_5A^3_7$, отрезками $P_5Q_5$ и $Q_5Q_6$ и лучом $Q_6A^0_3 = A^3_1A^0_3$, а $\alpha_6$ есть угол между лучами $A^3_1A^3_0 = A^4_7A^3_0$ и $A^3_1A^4_8 = A^4_7A^4_8$. Другими словами, пусть $\alpha_1$, $\alpha_2$, \ldots, $\alpha_6$ есть фигуры, на которые лучи $A_3A_2$, $A_4A_3$, \ldots, $A_7A_6$ разбивают угол $P_2P_1Q_2$. Пусть также точки $O_1$, $O_2$, \ldots, $O_5$ суть точки пересечения биссектрисы угла $P_2P_1Q_2$ с биссектрисами углов $P_1A_2Q_2$, $P_2A_3Q_3$, \ldots, $P_5A_6Q_6$ соответственно.

\begin{Lm}
Пусть точка $p$ лежит в углу $P_2P_1Q_2$. Тогда

$p \in \{O_1, O_2, \ldots, O_5\} \Leftrightarrow ((T'(p) \text{ определено}) \wedge (T'(p) = p))$.
\end{Lm}

\begin{Lm}
Индуцированное преобразование $T'$ есть кусочно-аффинное преобразование угла $P_2P_1Q_2$, устроенное таким образом, что:
\begin{itemize}
    \item для фигуры $\alpha_i$, $i \in \{1,2,3,4,5\}$, $T'$ есть поворот на угол $\frac{(6-i)\pi}{6}$ против часовой стрелки вокруг точки $O_i$;
    \item для фигуры-угла $\alpha_6$, $T'$ есть параллельный перенос на вектор $\overrightarrow{A^3_1A^3_7}$.
\end{itemize}

\end{Lm}

Для удобства будем считать, что на границах фигур $\alpha_i$, $i \in \{1,2,\ldots,6\}$, преобразование $T'$ не определено.

\begin{Lm}
Пусть точка $p \in Z'$. Тогда $p$ является периодической (граничной, апериодической) относительно преобразования внешнего биллиарда $T$, если и только если $p$ есть периодическая (граничная, апериодическая) точка относительно преобразования $T'$.
\end{Lm}

Введем также <<индуцированный>> код $\rho'$ для преобразования $T'$.

\begin{Def}
Пусть $p \in Z'$ - периодическая или апериодическая точка. Тогда $\rho'(p)$ есть бесконечная в обе стороны последовательность $(u'_n), n \in \mathbb{Z}$, т.ч. $\forall n \in \mathbb{Z}: T'^{n}(p) \in int(\alpha_{u'_n})$.
\end{Def}

Аналогичным образом можно определить и <<индуцированный>> код и для граничных точек. Далее под словом <<код>> мы будем подразумевать именно $\rho'$, если не указано обратное.

Введем также еще одно техническое определение.

\begin{Def}
Пусть $p \in Z'$ - точка с кодом $(u'_n)$, $n \in \mathbb{Z} \cap (l_1, l_2)$, $-\infty \leq l_1 < l_2 \leq +\infty$. Тогда будем обозначать $u'_n$, $n \in (l_1, l_2)$, как $\rho(p)[n]$, а подпоследовательность $u'_{k_1}u'_{k_1+1}\ldots u'_{k_2}$, $l_1 < k_1 \leq k_2 < l_2$ как $\rho'(p)[k_1..k_2] \equiv \rho'(p)[k_1, k_2]$.
\end{Def}

\section{Поиск периодических компонент}

Введем понятие периодической компоненты.

\begin{Def}
Компонентой называется максимальное по включению связное множество точек с одинаковым кодом. Компонента называется периодической, если хотя бы одна точка в ней является периодической.
\end{Def}

Сформулируем несколько лемм об устройстве периодических компонент.

\begin{Lm}
Точка $p$ является периодической, если и только если периодическим является ее код.
\end{Lm}

Периодические компоненты устроены следующим образом.

\begin{Lm} \label{lm:periodic1}
Периодическая компонента есть открытый выпуклый невырожденный многоугольник, стороны которого параллельны сторонам стола $\gamma$.
Граница же периодической компоненты состоит исключительно из граничных точек.
\end{Lm}

Введем также понятие периода.

\begin{Def}
Периодом точки $p$ $per_T(p)$ / компоненты $U$ $per_T(U)$ относительно преобразования $T$ назовем минимальное натуральное число $n$, т.ч. $T^n(p) = p$ / $T^n(U) = U$, если такое $n$ существует. Аналогичным образом введем и период $per_{T'}$ относительно преобразования $T'$.
\end{Def}

Следующая лемма устанавливает важные свойства периодических компонент.

\begin{Lm} \label{lm:periodic2}
    Для любой периодической компоненты $U$ выполнены следующие утверждения:
\begin{enumerate}
    \item Все точки $U$ периодические.
    \item $U$ имеет период, как $per_T$, так и $per_{T'}$.
    \item Если $U$ есть центрально-симметричный многоугольник c центром $c$, и $per_T(c)$ нечетен, то $per_T(c) = per_T(U)$, и $\forall p \in U \backslash \{c\}: per_T(p) = 2per_T(U)$; в любом ином случае, $\forall p \in U: per_T(p) = per_T(U)$.
    \item $T'^{per_{T'}(U)}(U)$ есть поворот $U$ на угол $\frac{\pi l}{6}$ для некоторого $l \in \mathbb{Z}$ вокруг точки $c \in U$, являющейся центром масс $U$, причем $per_{T'}(c) = per_{T'(U)}$, и $\forall p \in U \backslash \{c\}: per_{T'}(p) = per_{T'}(U) * 12 / \text{НОД(l, 12)}$.
\end{enumerate} 
\end{Lm}

Для любой заданной периодической точки $p$, лежащей внутри угла $P_2P_1Q_2$ с координатами в $\mathbb{Q}[\sqrt{3}]$ можно найти ее периодическую {\it компоненту} с помощью алгоритма \ref{alg:findPeriodicComponent}.

\begin{algorithm}
\caption{Поиск периодической компоненты, содержащей заданную точку.}
\label{alg:findPeriodicComponent}
\KwData {периодическая точка $startPoint$, лежащая в угле $P_2P_1Q_2$;}
\KwResult {многоугольник $U$, являющийся периодической компонентой, содержащей точку $p$.}

\SetKwProg{Fn}{Function}{;}{end function}
\Fn{findPeriodicComponent($startPoint$)}{
    $p := startPoint$
    $U_0 := \text{угол }P_2P_1Q_2$\;
    $p_0 := p$\;
    \tcc{бесконечный цикл}
    \For{$i = 0, 1, 2, \ldots$}{
        $p := T'(p)$\;
        $U :=$ той (возможно, бесконечно-)многоугольной части $T'(U)$, которая содержит $p$\;
        \If{многоугольник $U$ уже встречался ранее}{
            {\bf break}\;
        }
    }
    \While{точка $startPoint$ НЕ лежит внутри $U$}{
        $U := T'(U)$\;
    }
    \Return {$U$}\;
}

\end{algorithm}

Именно с помощью такого алгоритма, запущенного на ЭВМ, можно получить следующий результат. Будем говорить, что многоугольник $A$ вписан в многоугольник $B$, если на каждой из сторон $B$ целиком лежит хотя бы одна из сторон $A$.

\begin{Lm}
\begin{enumerate}

    \item Периодической компонентой с кодом, равным $\ldots 11111 \ldots$, является правильный двенадцатиугольник $W_1$, вписанный в треугольник $P_2P_1Q_2$, c центром в точке $O_1$.
    
    \item Периодической компонентой с кодом, равным $\ldots 22222 \ldots$, является равносторонний, но неправильный шестиугольник $W_2$ с углами $\frac{\pi}{2}, \frac{5\pi}{6}, \frac{\pi}{2}, \frac{5\pi}{6}, \frac{\pi}{2}, \frac{5\pi}{6}$ и центром $O_2$, вписанный в четырехугольник $P_3P_2Q_2Q_3$ таким образом, что вершина одного из прямых углов совпадает с точкой $P_3$, а противоположная вершина - с точкой $Q_2$.
    
    \item Периодической компонентой с кодом, равным $\ldots 33333 \ldots$, является равносторонний, но неправильный восьмиугольник $W_3$ с углами $\frac{2\pi}{3}, \frac{5\pi}{6}, \frac{2\pi}{3}, \frac{5\pi}{6}, \frac{2\pi}{3}, \frac{5\pi}{6}, \frac{2\pi}{3}, \frac{5\pi}{6}$ и центром $O_3$, вписанный в четырехугольник $P_4P_3Q_3Q_4$ таким образом, что вершина одного из углов, равных $\frac{2\pi}{3}$, совпадает с точкой $Q_3$.
    
    \item Периодической компонентой с кодом, равным $\ldots 44444 \ldots$, является правильный двенадцатиугольник $W_4$ с центром в $O_4$, вписанный в четырехугольник $P_5P_4Q_4Q_5$.
    
\end{enumerate}

\end{Lm}

\begin{figure}[h!]
\begin{center}
\includegraphics[width=100mm]{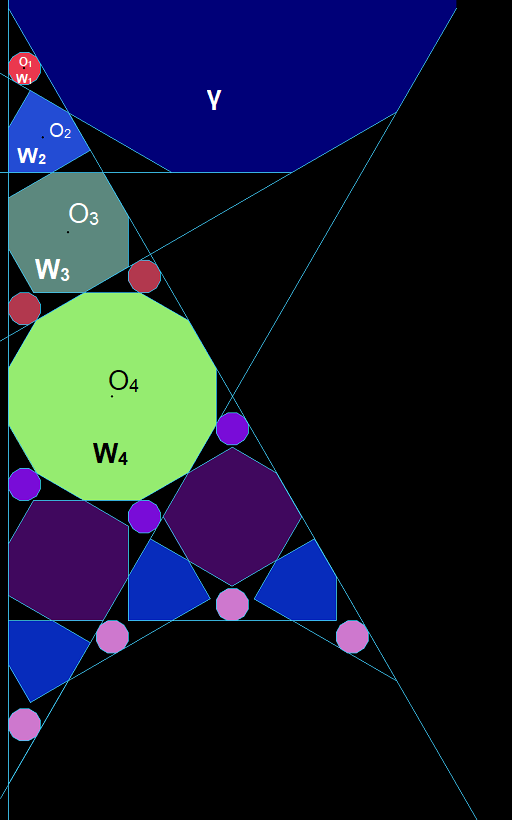}
\caption{Базовые периодические компоненты}
\label{pic:periodicComponents}
\end{center}
\end{figure}

Таким же методом можно найти и другие периодические компоненты; некоторые из них изображены на рис. \ref{pic:periodicComponents}. Отметим, что в случае пятиугольника \cite{Tab93}, все периодические компоненты (как минимум, внутри первой инвариантной компоненты) являются правильными пятиугольниками и десятиугольниками; в случае восьмиугольника \cite{Rukhovich18}, все периодические компоненты являются правильными восьмиугольниками. Случай же двенадцатиугольника оказывается сложнее, ибо уже среди <<базовых>> периодических компонент встречаются неправильные многоугольники, причем возможно, что это еще не все реально встречающиеся типы периодических компонент!

Тем не менее, пока мы вполне придерживаемся намеченного в \cite{Rukhovich18} плана, следующим шагом которого является поиск преобразования первого возвращения для некоторых фигур.

\section{Преобразования первого возвращения}

\begin{Def}
Пусть некоторая фигура $S$ лежит внутри угла $P_2P_1Q_2$. Тогда определим преобразование $T'_S$ первого возвращения (first return map) на $S$ относительно преобразования $T'$ таким образом, что $\forall s \in S$: $T'_S(s) := T'^{n_s}(s)$, где $n_s$ есть минимальное целое положительное число, т.ч. $T'^{n_s}(s) \in S$. Если такого $n_s$ не существует для некоторой точки $s \in S$, то для этой точки $T'_S(s)$ не определено.
\end{Def}

Мы будем искать преобразования первого возвращения для фигур $S$ специального вида.

\begin{Def}
Пусть фигура $S$ лежит внутри угла $P_2P_1Q_2$. Будем говорить, что $S$ имеет красивое самовозвращение, если выполнено следующее условие:

\begin{itemize}
\item Пусть $x, y$ - две неграничные точки, лежащие в $S$, а $k$ - некоторое целое положительное число. Пусть $T'^k(x)$ и $T'^k(y)$ определены, причем $\forall i \in [1, k)$: $T'^i(x), T'^i(y) \notin S$. Пусть также последовательности $\rho'(x)[0, k-1]$ и $\rho'(y)[0, k-1]$ совпадают. Тогда точки $T'^k(x)$ и $T'^k(y)$ либо обе лежат внутри $S$, либо обе не лежат внутри $S$.
\end{itemize}

\end{Def}

Таким свойством, в частности, обладает любой многоугольник, выпуклый или невыпуклый, который ограничен со всех сторон сторонами угла $P_2P_1Q_2$ и/или периодическими компонентами с траекториями, не пересекающими (лишь касающимися внешним образом) самого исходного многоугольника. Для такого многоугольника, корректным методом поиска преобразования первого возвращения является алгоритм \ref{alg:findFirstReturnMap}.

\begin{algorithm}
\caption{Поиск преобразования первого возвращения}
\label{alg:findFirstReturnMap}
\KwData {открытый многоугольник $S$, имеющий красивое самовозвращение}
\KwResult {два массива открытых многоугольников $S_1 = S_1[0..l-1]$ и $S_2 = S_2[0..l-1]$, т.ч.:
\begin{enumerate}
    \item $\forall i \in [0, l)\ \exists $ движение $f_i$, т.ч. $f_i(S_1[i]) = S_2[i]$, причем для всех точек $p \in int(S_1[i]): T'_S(p) = f_i(p)$;
    \item $\forall i \in [0, l)$ границы $S_1[i]$ и $S_2[i]$ состоят исключительно из граничных точек.
\end{enumerate}
}

\SetKwProg{Fn}{Function}{;}{end function}
\Fn{tryFirstReturnMap($S$)}{
    $S_1 := S_2 := []$\;
    $pols = [S]$\;
    \Repeat{$pols$ пуст}{
        \ForEach{$pol \in pols$}{
            разбить $pol$ лучами $A_3A_2$, $A_4A_3$, \ldots, $A_7A_6$ на многоугольники и поместить их(многоугольники) в массив $newPols$\;
            $npols = []$\;
            \ForEach{$newPol \in newPols$}{
                $newPol = T'(newPol)$\;
                добавить $newPol$ в $S_2$, если $newPol \subset S$, и в $npols$ иначе\;
            }
        }
        $pols := npols$;
    }
    \ForEach{$pol \in S_2$}{
        последовательно применять $T'^{-1}$ к $pol$, пока $pol$ не окажется внутри $S$ (не менее одного раза); результат добавить в $S_1$\; 
    }
    \Return $(S_1, S_2)$\;
}

\end{algorithm}

Отметим, что теоретически, такой алгоритм может не завершиться никогда; однако если же он завершится, то при обладании входного многоугольника красивым самовозвращением, преобразование первого возвращения будет построено корректно, ибо каждый из промежуточных многоугольников будет либо целиком лежать внутри $S$, либо не пересекаться с ним.

\section{Самоподобие}

Согласно намеченному плану, рассмотрим преобразование $\Gamma_1$, являющееся сжатием с центром в точке $A_1 = P_1$ и переводящее $O_5$ в $O_1$. Пусть $Z'_1 := \Gamma_1(Z')$. По аналогии с \cite{Rukhovich18}, ожидается, что $T'_{Z'_1}$ будет похож на $T'$ для $Z'$. Однако алгоритм \ref{alg:findFirstReturnMap} сообщает, что верна

\begin{Lm}
$T'_{Z'_1}$ разбивает $Z'_1$ на десять многоугольников, четыре из которых суть треугольники, пять - четырехугольники, а один - неравносторонний шестиугольник.
\end{Lm}

Таким образом, план в чистом виде реализовать не удалось. Однако самоподобие найти все-таки удается, хоть и более сложным способом. Для этого введем сжатие $\Gamma_3$ с центром в той же $A_1$, но переводящее $O_5$ в $O_4$; пусть $Z'_4 = \Gamma_4(Z')$, а $Z'_{14} = \Gamma_1(Z'_4)$. Заметим (с помощью компьютера), что $\Gamma_1(W_4)$ есть периодическая компонента с периодом 37 относительно $T'$. Следовательно, $Z'_4$ и $Z'_{14}$, как и $Z'_1$ обладают красивым самовозвращением, и к ним можно применить алгоритм \ref{alg:findFirstReturnMap}.

\begin{Lm} \label{lm:similarFRMs}
Пусть $S = Z'_4$ или $S = Z'_{14}$. Тогда $T'_S$ разбивает $S$ на восемь многоугольников, два из которых - треугольники, а остальные шесть - четырехугольники, причем один из этих четырехугольников невыпуклый. Более того, разбиение $Z'_{14}$ можно получить из разбиения $Z'_4$, применив к многоугольникам преобразование $\Gamma_1$.
\end{Lm}

Прямым следствием леммы \ref{lm:similarFRMs} является

\begin{Lm}
    Пусть $p_4 \in int(Z'_4)$, a $p_{14} := \Gamma_1(p_4)$. Тогда:
    \begin{enumerate}
        \item $T'_{Z'_4}(p_4)$ определено, если и только если $T'_{Z'_{14}}(p_{14})$ определено;
        \item если $T'_{Z'_4}(p_4)$ определено, то $\Gamma_1(T'_{Z'_4}(p_4)) = T'_{Z'_{14}}(p_{14}) = T'_{Z'_{14}}(\Gamma_1(p_4))$.
    \end{enumerate}
\end{Lm}

\section{Доказательство теоремы \ref{MainTreorem}}

\begin{figure}[h!]
\begin{center}
\includegraphics[width=100mm]{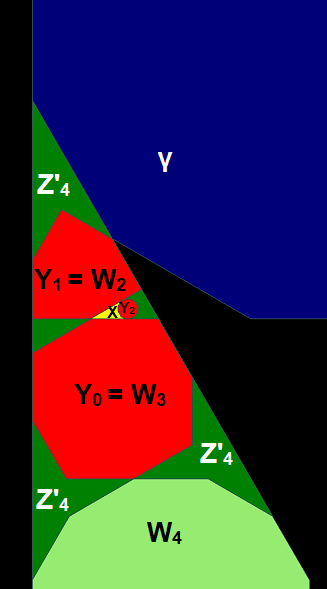}
\caption{Фигура $X$ и начало последовательности $Y$}.
\label{pic:spiral}
\end{center}
\end{figure}

 C помощью найденного в прошлом разделе самоподобия, докажем существование апериодической траектории, т.е. теорему \ref{MainTreorem}. Для этого рассмотрим (см. рис. \ref{pic:spiral}) <<ракету>> $X := (T')^{-2}(Z'_{14})$. Отметим, что $X$ ограничена фигурами $W_3$, $W_2$ и $(T')^{-2}(\Gamma_1(W_4))$, т.е. тремя периодическими фигурами. Пусть $\Gamma_X$ есть аффинное преобразование, переводящее $Z_4$ в $X$ с сохранением ориентации. Так как $(T')^{-2}$ не разделило <<ракету>> $Z'_{14}$, то очевидна следующая, аналогичная предыдущей, лемма.

\begin{Lm} \label{lm:selfSimilarity}
    Пусть $p_4 \in int(Z'_4)$, a $p_X := \Gamma_X(p_4)$. Тогда:
    \begin{enumerate}
        \item $T'_{Z'_4}(p_4)$ определено, если и только если $T'_X(p_X)$ определено;
        \item если $T'_{Z'_4}(p_4)$ определено, то $\Gamma_X(T'_{Z'_4}(p_4)) = T'_X(p_X) = T'_X(\Gamma_X(p_4))$.
    \end{enumerate}
\end{Lm}

Рассмотрим бесконечную последовательность фигур ($Y_0$, $Y_1$, $Y_2$, \ldots), где:

\begin{itemize}
    \item $Y_0 = W_3$, $Y_1 = W_2$, $Y_2$ = $(T')^{-2}(\Gamma_1(W_4))$;
    \item $\forall n \in \mathbb{N}, n \geq 3: Y_n = \Gamma_X(Y_{n-3})$.
\end{itemize}

Легко показать, что все фигуры $Y_n$ являются периодическими, причем так как $T'_{Z'_4}(X)$ не пересекается с $X$, то относительно $T'_{Z'_4}$, период $Y_n$ как минимум в два раза больше, нежели период $Y_{n-3}$, для любого $n \geq 3$. Пусть $y$ есть предельная точка последовательности ($Y_n$). Тогда $y$ не может быть периодической точкой, ибо тогда по леммам \ref{lm:periodic1}, \ref{lm:periodic2} $y$ должна обладать окрестностью точек с одинаковыми периодами (кроме, быть может, одной). Но $y$ также не может быть и граничной точкой, ибо тогда должен существовать отрезок прямой ненулевой длины, содержащий $y$ и состоящий лишь из граничных точек; это невозможно в силу <<спиралеобразности>> последовательности $(Y_n)$. Следовательно, $y$ есть апериодическая точка, и теорема \ref{MainTreorem} доказана.

\section{Доказательства теорем \ref{MainTreorem2} и \ref{MainTreorem3} }

Доказательства теорем \ref{MainTreorem2} и \ref{MainTreorem3}, базируются на следующих леммах.

\begin{Lm} \label{lm:periodicComponentsOfRank0}
Пусть $S = Z'_4$ или $S = Z'_{14}$. По лемме \ref{lm:similarFRMs} преобразование первого возвращения $T'_S$ разбивает $S$ на восемь фигур; пусть это $S_1$, $S_2$, \ldots, $S_8$, а $t_1$, $t_2$, \ldots, $t_8$ - их <<времена>> первого возвращения в $S$. Тогда в $Z'$ существует $n$ периодических фигур $R_1$, \ldots, $R_n$, с периодами $p_1$, $p_2$, \ldots, $p_n$ относительно $T'$, т.ч. существует разбиение

$Z' = (\bigcup\limits_{i = 1}^{8} \bigcup\limits_{j = 0}^{t_i - 1} T'^{j}(S_i))) \cup (\bigcup\limits_{i = 1}^{n} \bigcup\limits_{j = 0}^{p_i - 1} T'^{j}(R_i)))$.

При этом, $n = 7$ при $S = Z'_4$ и $n = 20$ при $S = Z'_{14}$.   
\end{Lm}

\begin{Lm} \label{lm:figuresOfMeasure0}
Пусть $S$ есть произвольный многоугольник. Рассмотрим последовательность разбиений $S$ на конечное число многоугольников, устроенную следующим образом:

\begin{enumerate}
\item каждый многоугольник покрашен в зеленый или красный цвет;
\item каждое следующее разбиение может быть получено из предыдущего путем разбиения некоторых {\bf зеленых} многоугольников на некоторое количество красных и зеленых многоугольников, без изменения красных многоугольников;
\item существуют вещественное число $\epsilon > 0$ и целое положительное число $k$ такое, что для любого целого положительного $i$ и любого зеленого многоугольника $U$, участвующего в $i$-ом разбиении, гарантируется, что в $i+k$-ом разбиении красные фигуры, лежащие внутри $U$, обладают суммарной площадью не меньшей, чем $\epsilon A$, где $A$ есть площадь фигуры $U$.
\end{enumerate}

Тогда участвующие в разбиениях красные фигуры образуют в $S$ множество полной меры.
\end{Lm}

Утверждение леммы \ref{lm:periodicComponentsOfRank0} может быть проверено с помощью алгоритма, похожего на алгоритм \ref{alg:findFirstReturnMap}. В ходе проверки автором были обнаружены периодические компоненты с периодами 1, 1, 18, 24, 1, 60, 54, 3, 32, 2, 756, 1008, 48, 1, 2, 3, 4, 37, 42, 85 (периоды относительно $T'$).

Леммы \ref{lm:selfSimilarity}, \ref{lm:periodicComponentsOfRank0} и \ref{lm:figuresOfMeasure0} позволяют доказать теоремы \ref{MainTreorem2} и \ref{MainTreorem3} для фигуры $Z'$, ибо:

\begin{itemize}
    \item $Z'$ можно разбить на траектории первого возвращения точек фигуры $Z'_4$ и периодические компоненты, траектории которых на попадают в $Z'_4$ (лемма \ref{lm:periodicComponentsOfRank0});
    \item $Z'_4$, как и $Z'$, также можно разбить на траектории первого возвращения точек фигуры $S'_{14}$ и периодические компоненты, траектории которых не попадают в $Z'_{14}$ (та же лемма \ref{lm:periodicComponentsOfRank0});
    \item для любой периодической компоненты $U \subset Z'_{4}$ существует периодическая компонента $U_b \subset (Z'_{4} \ Z'_{14})$ с траекторией, не проходящей через $Z'_{14}$, а также целые неотрицательные числа $i, j$, т.ч. $U = T'^j(\Gamma_1^i(U_b))$ (лемма \ref{lm:selfSimilarity});
    \item можно рассмотреть последовательность разбиений на траектории первого возвращения фигуры $\Gamma_1^n(Z_4)$ и периодические компоненты, $n = 1, 2, \ldots$ (лемма \ref{lm:selfSimilarity}); если в каждом разбиении красить периодические компоненты в красный цвет, а траектории первого возвращения - в зеленый, то к полученной последовательности разбиений можно применить лемму \ref{lm:figuresOfMeasure0} с $k = 2$ (в отличие от случаев правильных восьмиугольника и пятиугольника, где $k = 1$);
    \item если ввести код $\rho''$ для точек $Z'_4$ и преобразования $T'_{Z'_4}$, кодирующий каждую точку числом от 1 до 8, то для любой периодической точки (компоненты) $p$, $\rho''(\Gamma(p))$ можно получить из $\rho''(p)$ с помощью подстановки $\sigma$, не зависящей от выбора точки $p$. 
\end{itemize}

Доказать же теоремы \ref{MainTreorem2} и \ref{MainTreorem3} для всего угла $P_2P_1Q_2$ и, как следствие, для всей плоскости, позволяет следующая лемма, похожая на лемму 10 в \cite{Rukhovich18}. Пусть $T_6$ есть преобразование первого возвращения $T'$ для угла $\alpha_6$. Пусть преобразование $H$ есть параллельный перенос, т.ч. $H(A_1) = A^3_1$.  

\begin{Lm} \label{selfSimilarity4}
Пусть $x \in int(\angle P_2P_1Q_2)$. Тогда:
\begin{enumerate}
\item $T'(x)$ определено, если и только если $T_6(H(x))$ определено;
\item Если $T'(x)$ определено, то $H(T'(x)) = T_6(H(x))$. 
\end{enumerate}
\end{Lm}

О результатах данной статьи были сделаны доклады на Combinatorics on Words, Calculability and Automata research school (CIRM, Marseille, France, January 30 - 3 February 2017), на Tiling Dynamical System research school (CIRM, Marseille, France, 20-24 November 2017), на XXIV и XXV Международных конференциях с международным участием "Ломоносов" (10-14 апреля 2017, 3-9 апреля 2018), на 57, 58 и 59 конференциях с международным участием в МФТИ (24-29 ноября 2014, 23-28 ноября 2015, 21-26 ноября 2016), а также на Зимней школе "Комбинаторика и теория алгоритмов" (18-25 февраля 2018).

Работа поддержана грантом РНФ № 17-11-01337.

\section{Литература}

\bibliographystyle{utf8gost705u}
\bibliography{biblio}

\end{document}